\documentclass[smallextended,numbook,runningheads]{svjour31}     
\smartqed  
\usepackage{graphicx}
\usepackage{mathptmx}      

\usepackage{latexsym,bm,amssymb,amsfonts,amsbsy,amsmath,graphics,color}

\newtheorem{theo}{Theorem}
\newtheorem{lem}{Lemma}

\newtheorem{rem}{Remark}

\newtheorem{defn}{Definition}
\def\RR{\mathbb R}
\def\CC{\mathbb C}
\def\pmatrix{ \left( \begin{array} }
\def\endpmatrix{ \end{array} \right) }

\def\bfy{{\bf y}}

\def\bfo{{\bf 0}}

\def\no{\noindent}
\def\diag{{\rm diag}}
\def\proof{\underline{Proof}\quad}
\def\QED{~\mbox{$\Box$}}

\def\phi{\varphi}

\def\P{{\cal P}}

\def\I{{\cal I}}
\def\P{{\cal P}}

\def\sigmd{{\dot\sigma}}
\def\O{\Omega}

\def\dd{\mathrm{d}}

\journalname{~}
\begin{document}

\title{Isospectral Property of Hamiltonian Boundary
Value Methods (HBVMs) and their connections with Runge-Kutta
collocation methods\thanks{Work developed within the project
``Numerical methods and software for
differential equations''.}} 

\titlerunning{HBVMs and collocation methods}     

\author{L.\,Brugnano \and
F.\,Iavernaro \and
 D.\,Trigiante }

\authorrunning{B.I.T.} 

\institute{Luigi Brugnano \at Dipartimento di Matematica,
Universit\`a di Firenze, Viale Morgagni 67/A, 50134 Firenze
(Italy).\\ \email{luigi.brugnano@unifi.it}
\and
Felice Iavernaro \at Dipartimento di Matematica, Universit\`a di Bari,
Via Orabona  4,  70125 Bari (Italy).\\ \email{
felix@dm.uniba.it}
\and
Donato Trigiante \at Dipartimento di Energetica, Universit\`a di Firenze,
Via Lombroso 6/17, 50134 Firenze (Italy).\\ \email{
trigiant@unifi.it} }

\date{Date: February 23, 2010.}

\maketitle

\begin{abstract}
One main issue, when numerically integrating autonomous Hamiltonian
systems, is the long-term conservation of some of its invariants,
among which the Hamiltonian function itself. Recently, a new class
of methods, named {\em Hamiltonian Boundary Value Methods (HBVMs)}
has been introduced and analysed \cite{BIT}, which are able to
exactly preserve polynomial Hamiltonians of arbitrarily high degree.
We here study a further property of such methods, namely that of
having, when cast as a Runge-Kutta method, a matrix of the Butcher
tableau with the same spectrum (apart from the zero eigenvalues) as
that of the corresponding Gauss-Legendre method, independently of
the considered abscissae. Consequently, HBVMs are always perfectly
$A$-stable methods. This, in turn, allows to elucidate the existing
connections with classical Runge-Kutta collocation methods.

\end{abstract}

\section{Introduction}\label{intro}

Hamiltonian problems are of great interest in many fields of
application, ranging from the macro-scale of celestial mechanics, to
the micro-scale of molecular dynamics. They have been deeply
studied, from the point of view of the mathematical analysis, since
two centuries. Their numerical solution is a more recent field of
investigation, which has led to define symplectic methods, i.e., the
simplecticity of the discrete map, considering that, for the
continuous flow, simplecticity implies the conservation of $H(y)$.
However, the conservation of the Hamiltonian and simplecticity of
the flow cannot be satisfied at the same time unless the integrator
produces the exact solution (see \cite[page\,379]{HLW}). More
recently, the conservation of energy has been approached by means of
the concept of the {\em discrete line integral}, in a series of
papers \cite{IP1,IP2,IT1,IT2,IT3}, leading to the definition of {\em
Hamiltonian Boundary Value Methods (HBVMs)}
\cite{BIS,BIT2,BIT,BIT1}, which is a class of methods able to
preserve, for the discrete solution, polynomial Hamiltonians of
arbitrarily high degree (and then, a {\em practical} conservation of
any sufficiently differentiable Hamiltonian). In more details, in
\cite{BIT}, HBVMs based on Lobatto nodes have been analysed, whereas
in \cite{BIT1} HBVMs based on Gauss-Legendre abscissae have been
considered. In the last reference, it has been actually shown that
both formulae are essentially equivalent to each other, since the
order and stability properties of the method turn out  to be
independent of the abscissae distribution, and both methods are
equivalent, when the number of the so called {\em silent stages}
tends to infinity. In this paper this conclusion if further
supported, since we prove that HBVMs, when cast as Runge-Kutta
methods, are such that  the corresponding matrix of the tableau has
the nonzero eigenvalues coincident with those of the corresponding
Gauss-Legendre formula (isospectral property of HBVMs).

This property can be also used to further analyse the existing
connections between HBVMs and Runge-Kutta collocation methods.

With this premise, the structure of the paper is the following: in
Section~\ref{hbvms} the basic facts about HBVMs are recalled; in
Section~\ref{iso} we state the main result of this paper,
concerning the isospectral property; in Section~\ref{coll} such
property is further generalized to study the existing connections
between HBVMs and Runge-Kutta collocation methods; finally, in
Section~\ref{fine} a few concluding remarks are given.

\section{Hamiltonian Boundary Value Methods}\label{hbvms}

The arguments in this section are worked out starting from the
arguments used in \cite{BIT,BIT1} to introduce and analyse HBVMs.
We consider canonical Hamiltonian problems in the form
\begin{equation}\label{hamilode}\dot y = J\nabla H(y),  \qquad y(t_0) = y_0\in\RR^{2m},
\end{equation}

\no where $J$ is a skew-symmetric constant matrix, and the
Hamiltonian $H(y)$ is assumed to be sufficiently differentiable.
The key formula which HBVMs rely on, is the {\em line integral}
and the related property of conservative vector fields:
\begin{equation}\label{Hy}
H(y_1) - H(y_0) = h\int_0^1 \sigmd(t_0+\tau h)^T\nabla
H(\sigma(t_0+\tau h))\dd\tau,
\end{equation}

\no for any $y_1 \in \RR^{2m}$, where $\sigma$ is any smooth
function such that
\begin{equation}
\label{sigma}\sigma(t_0) = y_0, \qquad\sigma(t_0+h) = y_1.
\end{equation}

\no Here we consider the case where $\sigma(t)$ is a polynomial of
degree $s$, yielding an approximation to the true solution $y(t)$
in the time interval $[t_0,t_0+h]$. The numerical approximation
for the subsequent time-step, $y_1$, is then defined by
(\ref{sigma}).

After introducing a set of $s$ distinct abscissae,
\begin{equation}\label{ci}0<c_{1},\ldots ,c_{s}\le1,\end{equation}

\no we set
\begin{equation}\label{Yi}Y_i=\sigma(t_0+c_i h), \qquad
i=1,\dots,s,\end{equation}

\no so that $\sigma(t)$ may be thought of as an interpolation
polynomial, interpolating the {\em fundamental stages} $Y_i$,
$i=1,\dots,s$, at the abscissae (\ref{ci}). We observe that, due
to (\ref{sigma}), $\sigma(t)$ also interpolates the initial
condition $y_0$.

\begin{rem}\label{c0} Sometimes, the interpolation at $t_0$ is
explicitly  required. In such a case, the extra abscissa $c_0=0$ is
formally added to (\ref{ci}). This is the case, for example, of a
Lobatto distribution of the abscissae \cite{BIT}.\end{rem}

Let us consider the following expansions of $\dot \sigma(t)$ and
$\sigma(t)$ for $t\in [t_0,t_0+h]$:
\begin{equation}
\label{expan} \dot \sigma(t_0+\tau h) = \sum_{j=1}^{s} \gamma_j
P_j(\tau), \qquad \sigma(t_0+\tau h) = y_0 + h\sum_{j=1}^{s}
\gamma_j \int_{0}^\tau P_j(x)\,\dd x,
\end{equation}

\no where $\{P_j(t)\}$ is a suitable basis of the vector space of
polynomials of degree at most $s-1$ and  the  (vector) coefficients
$\{\gamma_j\}$ are to be determined.  We shall consider an {\em
orthonormal basis} of polynomials on the interval
$[0,1]$\footnote{The use of an arbitrary polynomial basis is also
permitted and has been considered in the past (see for example
\cite{IT3,IP2}), however as was shown in \cite{BIT1}, among all
possible choices, the Legendre basis turns out to be the optimal
one.}, i.e.:
\begin{equation}\label{orto}\int_0^1 P_i(t)P_j(t)\dd t = \delta_{ij}, \qquad
i,j=1,\dots,s,\end{equation}

\no where $\delta_{ij}$ is the Kronecker symbol, and $P_i(t)$ has
degree $i-1$. Such a basis can be readily obtained as
\begin{equation}\label{orto1}P_i(t) = \sqrt{2i-1}\,\hat P_{i-1}(t), \qquad
i=1,\dots,s,\end{equation}

\no with $\hat P_{i-1}(t)$ the shifted Legendre polynomial, of
degree $i-1$, on the interval $[0,1]$.

\begin{rem}\label{recur}
From the properties of shifted Legendre polynomials (see, e.g.,
\cite{AS} or the Appendix in \cite{BIT}), one readily obtains that
the polynomials $\{P_j(t)\}$ satisfy the three-terms recurrence
relation:
\begin{eqnarray*}
P_1(t)&\equiv& 1, \qquad P_2(t) = \sqrt{3}(2t-1),\\
P_{j+2}(t) &=& (2t-1)\frac{2j+1}{j+1} \sqrt{\frac{2j+3}{2j+1}}
P_{j+1}(t) -\frac{j}{j+1}\sqrt{\frac{2j+3}{2j-1}} P_j(t), \quad
j\ge1.
\end{eqnarray*}
\end{rem}

We shall also assume that $H(y)$ is a polynomial, which implies
that the integrand in \eqref{Hy} is also a polynomial so that the
line integral can be exactly computed by means of a suitable
quadrature formula. It is easy to observe that in general, due to
the high degree of the integrand function,  such quadrature
formula cannot be solely based upon the available abscissae
$\{c_i\}$: one needs to introduce an additional set of abscissae
$\{\hat c_1, \dots,\hat c_r\}$, distinct from the nodes $\{c_i\}$,
in order to make the quadrature formula exact:
\begin{eqnarray} \label{discr_lin}
\displaystyle \lefteqn{\int_0^1 \sigmd(t_0+\tau h)^T\nabla
H(\sigma(t_0+\tau h))\mathrm{d}\tau   =}\\ && \sum_{i=1}^s \beta_i
\sigmd(t_0+c_i h)^T\nabla H(\sigma(t_0+c_i h)) + \sum_{i=1}^r \hat
\beta_i \sigmd(t_0+\hat c_i h)^T\nabla H(\sigma(t_0+\hat c_i h)),
\nonumber
\end{eqnarray}

\no where $\beta_i$, $i=1,\dots,s$, and $\hat \beta_i$,
$i=1,\dots,r$, denote the weights of the quadrature formula
corresponding to the abscissae $\{c_i\}\cup\{\hat c_i\}$, i.e.,
\begin{eqnarray}\nonumber
\beta_i &=& \int_0^1\left(\prod_{ j=1,j\ne i}^s
\frac{t-c_j}{c_i-c_j}\right)\left(\prod_{j=1}^r
\frac{t-\hat c_j}{c_i-\hat c_j}\right)\mathrm{d}t, \qquad i = 1,\dots,s,\\
\label{betai}\\ \nonumber \hat\beta_i &=& \int_0^1\left(\prod_{
j=1}^s \frac{t-c_j}{\hat c_i-c_j}\right)\left(\prod_{ j=1,j\ne
i}^r \frac{t-\hat c_j}{\hat c_i-\hat c_j}\right)\mathrm{d}t,
\qquad i = 1,\dots,r.
\end{eqnarray}

\begin{rem}\label{c01}
In the case considered in the previous Remark~\ref{c0}, i.e. when
$c_0=0$ is formally added to the abscissae (\ref{ci}), the first product in each
formula in (\ref{betai}) ranges from $j=0$ to $s$. Moreover, also the range
of $\{\beta_i\}$ becomes $i=0,1,\dots,s$. However, for sake of
simplicity, we shall not consider this case further.
\end{rem}

\begin{defn}\label{defhbvmks}
The method defined by the polynomial $\sigma(t)$, determined by
substituting the quantities in \eqref{expan} into the right-hand
side of \eqref{discr_lin}, and by choosing the unknown coefficient
$\{\gamma_j\}$ in order that the resulting expression vanishes, is
called {\em Hamiltonian Boundary Value Method with $k$ steps and
degree $s$}, in short {\em HBVM($k$,$s$)}, where $k=s+r$ \,
\cite{BIT}.\end{defn}

According to \cite{IT2}, the right-hand side of \eqref{discr_lin} is
called \textit{discrete line integral} associated with the map
defined by the HBVM($k$,$s$), while the vectors
\begin{equation}\label{hYi}
\hat Y_i \equiv \sigma(t_0+\hat c_i h), \qquad i=1,\dots,r,
\end{equation}

\no are called \textit{silent stages}: they just serve to increase,
as much as one likes, the degree of precision of the quadrature
formula, but they are not to be regarded as unknowns since, from
\eqref{expan} and (\ref{hYi}), they can be expressed in terms of
linear combinations of the fundamental stages (\ref{Yi}).

In the sequel, we shall see that HBVMs may be expressed through different,
though equivalent, formulations: some of them can be directly implemented in a
computer program, the others being of more theoretical interest.

Because of the equality \eqref{discr_lin}, we can apply the
procedure described in Definition \ref{defhbvmks} directly to the
original line integral appearing in the left-hand side. With this
premise, by considering the first expansion in \eqref{expan},  the
conservation property reads
$$\sum_{j=1}^{s} \gamma_j^T \int_0^1  P_j(\tau) \nabla
H(\sigma(t_0+\tau h))\dd\tau=0,$$

\no which, as is easily checked,  is certainly satisfied if we
impose the following set of orthogonality conditions:
\begin{equation}
\label{orth} \gamma_j = \int_0^1  P_j(\tau) J \nabla
H(\sigma(t_0+\tau h))\dd\tau, \qquad j=1,\dots,s.
\end{equation}

\no Then, from the second relation of \eqref{expan} we obtain, by
introducing the operator
\begin{eqnarray}\label{Lf}\lefteqn{L(f;h)\sigma(t_0+ch) =}\\
\nonumber && \sigma(t_0)+h\sum_{j=1}^s \int_0^c P_j(x) \dd x \,
\int_0^1 P_j(\tau)f(\sigma(t_0+\tau h))\dd\tau,\qquad
c\in[0,1],\end{eqnarray}

\no that $\sigma$ is the eigenfunction of $L(J\nabla H;h)$
relative to the eigenvalue $\lambda=1$:
\begin{equation}\label{L}\sigma = L(J\nabla H;h)\sigma.\end{equation}

\begin{defn} Equation (\ref{L}) is the {\em Master Functional
Equation} defining $\sigma$ ~\cite{BIT1}.\end{defn}

\begin{rem}\label{MFE}
From the previous arguments, one readily obtains that the Master
Functional Equation (\ref{L}) characterizes HBVM$(k,s)$ methods,
for all $k\ge s$. Indeed, such methods are uniquely defined by the
polynomial $\sigma$, of degree $s$, the number of steps $k$ being
only required to obtain an exact quadrature formula (see
(\ref{discr_lin})).\end{rem}

To practically compute $\sigma$, we set (see (\ref{Yi}) and
(\ref{expan}))
\begin{equation}
\label{y_i} Y_i=  \sigma(t_0+c_i h) = y_0+ h\sum_{j=1}^{s} a_{ij}
\gamma_j, \qquad i=1,\dots,s,
\end{equation}

\no where
$$a_{ij}=\int_{0}^{c_i} P_j(x) \mathrm{d}x, \qquad
i,j=1,\dots,s.$$%

\no Inserting \eqref{orth} into \eqref{y_i} yields the final
formulae which define the HBVMs class based upon the orthonormal
basis $\{P_j\}$:
\begin{equation}
\label{hbvm_int} Y_i=y_0+h  \sum_{j=1}^s a_{ij}\int_0^1
P_j(\tau)  J \nabla H(\sigma(t_0+\tau h))\,\dd\tau, \qquad
i=1,\dots,s.
\end{equation}

For sake of completeness, we report the nonlinear system
associated with the HBVM$(k,s)$ method, in terms of the
fundamental stages $\{Y_i\}$ and the silent stages $\{\hat Y_i\}$
(see (\ref{hYi})), by using the notation
\begin{equation}\label{fy}
f(y) = J \nabla H(y).
\end{equation}

\no In this context, it  represents the discrete counterpart of
\eqref{hbvm_int}, and may be directly retrieved by evaluating, for
example, the integrals in \eqref{hbvm_int} by means of the (exact)
quadrature formula introduced in \eqref{discr_lin}:
\begin{eqnarray}\label{hbvm_sys}
\lefteqn{ Y_i =}\\
&=& y_0+h\sum_{j=1}^s a_{ij}\left( \sum_{l=1}^s \beta_l
P_j(c_l)f(Y_l) + \sum_{l=1}^r\hat \beta_l P_j(\hat c_l) f(\widehat
Y_l) \right),\quad i=1,\dots,s.\nonumber
\end{eqnarray}

\no From the above discussion it is clear that, in the
non-polynomial case, supposing to choose the abscissae $\{\hat
c_i\}$ so that the sums in (\ref{hbvm_sys}) converge to an
integral as $r\equiv k-s\rightarrow\infty$, the resulting formula is
\eqref{hbvm_int}.

\begin{defn}\label{infh} Formula (\ref{hbvm_int}) is named {\em
$\infty$-HBVM of degree $s$} or {\em HBVM$(\infty,s)$} \, \cite{BIT1}.
\end{defn}

This implies that HBVMs may be as well applied in the non-polynomial case since,
in finite precision arithmetic, HBVMs are undistinguishable from their limit
formulae \eqref{hbvm_int}, when a sufficient number of silent stages is
introduced. The aspect of having a {\em practical} exact integral,
for $k$ large enough, was already stressed in
\cite{BIS,BIT,BIT1,IP1,IT2}.

On the other hand, we emphasize that, in the non-polynomial case,
\eqref{hbvm_int} becomes an {\em operative method} only after that a suitable
strategy to approximate the integrals appearing in it is taken into account. In
the present case, if one discretizes the {\em Master Functional Equation}
(\ref{Lf})--(\ref{L}), HBVM$(k,s)$ are then obtained, essentially
by extending the discrete problem (\ref{hbvm_sys}) also to the
silent stages (\ref{hYi}). In order to simplify the exposition, we
shall use (\ref{fy}) and introduce the following notation:
\begin{eqnarray*}
\{\tau_i\} = \{c_i\} \cup \{\hat{c}_i\}, &&
\{\omega_i\}=\{\beta_i\}\cup\{\hat\beta_i\},\\[2mm]
y_i = \sigma(t_0+\tau_ih), && f_i = f(\sigma(t_0+\tau_ih)), \qquad
i=1,\dots,k.
\end{eqnarray*}

\no The discrete problem defining the HBVM$(k,s)$ then becomes,
\begin{equation}\label{hbvmks}
y_i = y_0 + h\sum_{j=1}^s \int_0^{\tau_i} P_j(x)\dd x
\sum_{\ell=1}^k \omega_\ell P_j(\tau_\ell)f_\ell, \qquad
i=1,\dots,k.
\end{equation}

By introducing the vectors $$\bfy = (y_1^T,\dots,y_k^T)^T, \qquad
e=(1,\dots,1)^T\in\RR^k,$$ and the matrices
\begin{equation}\label{OIP}\O=\diag(\omega_1,\dots,\omega_k), \qquad
\I_s,~\P_s\in\RR^{k\times s},\end{equation} whose $(i,j)$th entry
are given by
\begin{equation}\label{IDPO}
(\I_s)_{ij} = \int_0^{\tau_i} P_j(x)\mathrm{d}x, \qquad
(\P_s)_{ij}=P_j(\tau_i), \end{equation}

\no we can cast the set of equations (\ref{hbvmks}) in vector form
as $$\bfy = e\otimes y_0 + h(\I_s
\P_s^T\O)\otimes I_{2m}\, f(\bfy),$$

\no with an obvious meaning of $f(\bfy)$. Consequently, the method
can be regarded as a Runge-Kutta method with the following Butcher
tableau:
\begin{equation}\label{rk}
\begin{array}{c|c}\begin{array}{c} \tau_1\\ \vdots\\ \tau_k\end{array} & \I_s \P_s^T\O\\
 \hline                    &\omega_1\, \dots~ \omega_k
                    \end{array}\end{equation}

\begin{rem}\label{ascisse} We observe that, because of the use of an
orthonormal basis, the role of the abscissae $\{c_i\}$ and of the
silent abscissae $\{\hat c_i\}$ is interchangeable, within the set
$\{\tau_i\}$. This is due to the fact that all the matrices
$\I_s$, $\P_s$, and $\O$ depend on all the abscissae $\{\tau_i\}$,
and not on a subset of them, and they are invariant with respect
to the choice of the fundamental abscissae $\{c_i\}$.
\end{rem}

In particular, when a Gauss distribution of the abscissae
$\{\tau_1,\dots,\tau_k\}$ is considered, it can be proved that the
resulting HBVM$(k,s)$ method \cite{BIT1}:

\begin{itemize}

\item has order $2s$ for all $k\ge s$;

\item is symmetric and perfectly $A$-stable (i.e., its
stability region coincides with the left-half complex plane,
$\CC^-$ \cite{BT});

\item reduces to the Gauss-Legendre method of order $2s$, when
$k=s$;

\item exactly preserves polynomial Hamiltonian functions of degree $\nu$,
provided that \begin{equation}\label{knu}k\ge \frac{\nu
s}2.\end{equation}

\end{itemize}

Additional results and references on HBVMs can be found at the
{\em HBVMs Homepage} \cite{HBVMsHome}.

\section{The Isospectral Property}\label{iso}

We are now going to prove a further additional result, related to
the matrix appearing in the Butcher tableau (\ref{rk}),
corresponding to HBVM$(k,s)$, i.e., the matrix
\begin{equation}\label{AMAT}A = \I_s \P_s^T\O\in\RR^{k\times
k}, \qquad k\ge s,\end{equation}

\no whose rank is $s$. Consequently it has a $(k-s)$-fold zero
eigenvalue. In this section, we are going to discuss the location
of the remaining $s$ eigenvalues of that matrix.

Before that, we state the following preliminary result, whose
proof can be found in \cite[Theorem\,5.6 on page\,83]{HW}.

\begin{lem}\label{gauss} The eigenvalues of the matrix
\begin{equation}\label{Xs}
X_s = \pmatrix{cccc}
\frac{1}2 & -\xi_1 &&\\
\xi_1     &0      &\ddots&\\
          &\ddots &\ddots    &-\xi_{s-1}\\
          &       &\xi_{s-1} &0\\
\endpmatrix, \end{equation} with
\begin{equation}\label{xij}\xi_j=\frac{1}{2\sqrt{(2j+1)(2j-1)}}, \qquad
j\ge1,\end{equation} coincide with those of the matrix in the
Butcher tableau of the Gauss-Legendre method of order
$2s$.\end{lem}

\medskip
We also need the following preliminary result, whose proof derives
from the properties of shifted-Legendre polynomials (see, e.g.,
\cite{AS} or the Appendix in \cite{BIT}).

\begin{lem}\label{intleg} With reference to the matrices in
(\ref{OIP})--(\ref{IDPO}), one has
$$\I_s = \P_{s+1}\hat{X}_s,$$

\no where
$$\hat{X}_s = \pmatrix{cccc}
\frac{1}2 & -\xi_1 &&\\
\xi_1     &0      &\ddots&\\
          &\ddots &\ddots    &-\xi_{s-1}\\
          &       &\xi_{s-1} &0\\
\hline &&&\xi_s\endpmatrix,$$

\no with the $\xi_j$ defined by (\ref{xij}). \end{lem}

\medskip
The following result then holds true.

\begin{theo}[Isospectral Property of HBVMs]\label{mainres}
For all $k\ge s$ and for any choice of the abscissae $\{\tau_i\}$
such that the quadrature defined by the weights $\{\omega_i\}$ is
exact for polynomials of degree $2s-1$, the nonzero eigenvalues of
the matrix $A$ in (\ref{AMAT}) coincide with those of the matrix
of the Gauss-Legendre method of order $2s$.
\end{theo}
\begin{proof}
For $k=s$, the abscissae $\{\tau_i\}$ have to be the $s$
Gauss-Legendre nodes, so that HBVM$(s,s)$ reduces to the Gauss
Legendre method of order $2s$, as outlined at the end of
Section~\ref{hbvms}.

When $k>s$, from the orthonormality of the basis, see
(\ref{orto}), and considering that the quadrature with weights
$\{\omega_i\}$ is exact for polynomials of degree (at least)
$2s-1$, one easily obtains that (see (\ref{OIP})--(\ref{IDPO}))
$$\P_s^T\O\P_{s+1} = \left( I_s ~ \bfo\right),$$

\no since, for all ~$i=1,\dots,s$,~ and ~$j=1,\dots,s+1$:
$$\left(\P_s^T\O\P_{s+1}\right)_{ij} = \sum_{\ell=1}^k \omega_\ell
P_i(t_\ell)P_j(t_\ell)=\int_0^1 P_i(t)P_j(t)\dd t = \delta_{ij}.$$

\no By taking into account the result of Lemma~\ref{intleg}, one
then obtains:
\begin{eqnarray}\nonumber
A\P_{s+1} &=& \I_s \P_s^T\O\P_{s+1} = \I_s \left(I_s~\bfo\right)
=\P_{s+1}
\hat{X}_s \left(I_s~\bfo\right) = \P_{s+1}\left(\hat{X}_s~\bfo\right)\\
 &=& \P_{s+1}
\pmatrix{cccc|c}
\frac{1}2 & -\xi_1 && &0\\
\xi_1     &0      &\ddots& &\vdots\\
          &\ddots &\ddots    &-\xi_{s-1}&\vdots\\
          &       &\xi_{s-1} &0&0\\
\hline &&&\xi_s&0\endpmatrix ~\equiv~ \P_{s+1}\widetilde X_s,
\label{tXs}
\end{eqnarray}

\no with the $\{\xi_j\}$ defined according to (\ref{xij}).
Consequently, one obtains that the columns of $\P_{s+1}$
constitute a basis of an invariant (right) subspace of matrix $A$,
so that the eigenvalues of $\widetilde X_s$ are eigenvalues of
$A$. In more detail, the eigenvalues of $\widetilde X_s$ are those
of $X_s$ (see (\ref{Xs})) and the zero eigenvalue. Then, also in
this case, the nonzero eigenvalues of $A$ coincide with those of
$X_s$, i.e., with the eigenvalues of the matrix defining the
Gauss-Legendre method of order $2s$.\QED\end{proof}

\section{HBVMs and Runge-Kutta collocation methods}\label{coll}

By using the previous results and notations, now we further
elucidate the existing connections between HBVMs and Runge-Kutta
collocation methods. We shall continue to use an orthonormal basis
$\{P_j\}$, along which the underlying {\em extended collocation}
polynomial $\sigma(t)$ is expanded, even though the arguments could
be generalized to more general bases, as sketched below. On the
other hand, the distribution of the internal abscissae can be
arbitrary.

Our starting point is a generic collocation method with $k$
stages, defined by the tableau
\begin{equation}
\label{collocation_rk}
\begin{array}{c|c}\begin{array}{c} \tau_1\\ \vdots\\ \tau_k\end{array} &  \mathcal A \\
 \hline                    &\omega_1\, \ldots  ~ \omega_k
\end{array}
\end{equation}

\no where, for $i,j=1,\dots,k$: $$\mathcal A=
(\alpha_{ij})\equiv\left(\int_0^{\tau_i} \ell_j(x) \mathrm{d}x
\right), \qquad \omega_j=\int_0^{1} \ell_j(x) \mathrm{d}x,$$

\no $\ell_j(\tau)$ being the $j$th Lagrange polynomial of degree
$k-1$ defined on the set of abscissae $\{\tau_i\}$.

Given a positive integer $s\le k$, we can consider a basis
$\{p_1(\tau), \dots, p_s(\tau)\}$ of the vector space of
polynomials of degree at most $s-1$, and we set
\begin{equation}
\label{P} \hat\P_s = \pmatrix{cccc}
p_1(\tau_1) & p_2(\tau_1) & \cdots & p_s(\tau_1) \\
p_1(\tau_2) & p_2(\tau_2) & \cdots & p_s(\tau_2) \\
\vdots   & \vdots   &        & \vdots \\
p_1(\tau_k) & p_2(\tau_k) & \cdots & p_s(\tau_k)
\endpmatrix_{k \times s}
\end{equation}

\no (note that $\hat\P_s$ is full rank since the nodes are
distinct). The class of Runge-Kutta methods we are interested in
is defined by the tableau
\begin{equation}
\label{hbvm_rk}
\begin{array}{c|c}\begin{array}{c} \tau_1\\ \vdots\\ \tau_k\end{array} &
A \equiv \mathcal A \hat\P_s \Lambda_s \hat\P_s^T \Omega\\
 \hline                    &\omega_1\, \ldots \ldots ~ \omega_k
\end{array}
\end{equation}

\no where $\Omega=\diag(\omega_1,\dots,\omega_k)$ (see
(\ref{OIP})) and $\Lambda_s=\diag(\eta_1,\dots,\eta_s)$; the
coefficients $\eta_j$, $j=1,\dots,s$, have to be selected by
imposing suitable consistency conditions on the stages $\{y_i\}$
(see, e.g., \cite{BIT1}). In particular, when the basis is
orthonormal, as we shall assume hereafter, then matrix $\hat\P_s$
reduces to matrix $\P_s$ in (\ref{OIP})--(\ref{IDPO}), $\Lambda_s
= I_s$, and consequently (\ref{hbvm_rk}) becomes
\begin{equation}
\label{hbvm_rk1}
\begin{array}{c|c}\begin{array}{c} \tau_1\\ \vdots\\ \tau_k\end{array} &
A \equiv \mathcal A \P_s \P_s^T \Omega\\
 \hline       &\omega_1\, \ldots \ldots ~ \omega_k
\end{array}
\end{equation}

We note that the Butcher array $A$ has rank which cannot exceed
$s$, because it is defined by {\em filtering} $\mathcal A$ by the
rank $s$ matrix $\P_s \P_s^T \Omega$.

The following result then holds true, which clarifies the existing
connections between classical Runge-Kutta collocation methods and
HBVMs.

\begin{theo}\label{collhbvm} Provided that the quadrature formula defined by the
weights $\{\omega_i\}$ is exact for polynomials at least $2s-1$
(i.e., the Runge-Kutta method defined by the tableau
(\ref{hbvm_rk1}) satisfies the usual simplifying assumption
$B(2s)$), then the tableau (\ref{hbvm_rk1}) defines a HBVM$(k,s)$
method based at the abscissae $\{\tau_i\}$.
\end{theo}

\proof Let us expand the basis $\{P_1(\tau),\dots,P_s(\tau)\}$
along the Lagrange basis $\{\ell_j(\tau)\}$, $j=1,\dots,k$,
defined over the nodes $\tau_i$, $i=1,\dots,k$: $$
P_j(\tau)=\sum_{r=1}^k P_j(\tau_r) \ell_r(\tau),
 \qquad j=1,\dots,s.$$

\no It follows that, for $i=1,\dots,k$ and $j=1,\dots,s$:
$$\int_0^{\tau_i} P_j(x) \mathrm{d}x = \sum_{r=1}^k P_j(\tau_r)
\int_0^{\tau_i} \ell_r(x) \mathrm{d}x = \sum_{r=1}^k P_j(\tau_r)
\alpha_{ir},$$

\no that is (see (\ref{OIP})--(\ref{IDPO}) and
(\ref{collocation_rk})),
\begin{equation}\label{APeqI}
\I_s = \mathcal A \P_s.
\end{equation}

\no By substituting (\ref{APeqI}) into (\ref{hbvm_rk1}), one
retrieves that tableau (\ref{rk}), which defines the method
HBVM$(k,s)$. This completes the proof.\QED

\medskip
The resulting Runge-Kutta method \eqref{hbvm_rk1} is then energy
conserving if applied to polynomial Hamiltonian systems
\eqref{hamilode} when the degree of $H(y)$, is lower than or equal
to a quantity, say $\nu$, depending on $k$ and $s$. As an example,
when a Gaussian distribution of the nodes $\{\tau_i\}$ is
considered, one obtains (\ref{knu}).

\begin{rem}[{\bf About Simplecticity}]\label{symplectic} The choice
of the abscissae $\{\tau_1,\dots,\tau_k\}$ at the Gaussian points
in $[0,1]$ has also another important consequence, since, in such
a case, the collocation method (\ref{collocation_rk}) is the Gauss
method of order $2k$ which, as is well known, is a {\em symplectic
method}. The result of Theorem~\ref{collhbvm} then states that,
for any $s\le k$, the HBVM$(k,s)$ method is related to the Gauss
method of order $2k$ by the relation: $$A = {\cal A}
(\P_s\P_s^T\O),$$

\no where the {\em filtering matrix} $(\P_s\P_s^T\O)$ essentially
makes the Gauss method of order $2k$ ``work'' in a suitable
subspace.
\end{rem}

It seems like the price paid to achieve such conservation property
consists in the lowering of the order of the new method with
respect to the original one \eqref{collocation_rk}. Actually this
is not true,  because a fair comparison would be to relate method
\eqref{rk}--\eqref{hbvm_rk1}  to a collocation method constructed
on $s$ rather than on $k$  stages, since the resulting nonlinear
system turns out to have dimension $s$, as shown in \cite{BIT}.
This computational aspect is fully elucidated in a companion paper
\cite{blend}, devoted to the efficient implementation of HBVMs,
where the Isospectral Property of the methods is fully exploited
for this purpose.

\subsection{An alternative proof for the order of HBVMs}

We conclude this section by observing that the order $2s$ of an
HBVM$(k,s)$ method, under the hypothesis that
\eqref{collocation_rk} satisfies the usual simplifying assumption
$B(2s)$, i.e., the quadrature defined by the weights
$\{\omega_i\}$ is exact for polynomials of degree at least $2s-1$,
may be stated by using an alternative, though equivalent,
procedure to that used in the proof of \cite[Corollary\,2]{BIT}
(see also \cite[Theorem\,2]{BIT1}).

Let us then define the $k \times k$ matrix $\P\equiv \P_k$  (see
(\ref{OIP})--(\ref{IDPO})) obtained by ``enlarging'' the matrix
$\P_s$ with $k-s$ columns defined by the normalized shifted
Legendre polynomials $P_j(\tau)$, $j=s+1,\dots,k$, evaluated at
$\{\tau_i\}$, i.e.,
$$\P=\pmatrix{ccc} P_1(\tau_1) & \dots &P_k(\tau_1)\\ \vdots &
&\vdots\\ P_1(\tau_k) & \dots &P_k(\tau_k)\endpmatrix.$$

\no By virtue of property $B(2s)$ for the quadrature formula
defined by the weights $\{\omega_i\}$, it satisfies
$$
\P^T \O \P =\pmatrix{ll} I_s & O \\ O & R
\endpmatrix, \qquad R\in\RR^{k-s\times k-s}.
$$

\no  This implies that $\P$ satisfies the property $T(s,s)$ in
\cite[Definition\,5.10 on page 86]{HW}, for the quadrature formula
$(\omega_i,\tau_i)_{i=1}^k$. Therefore, for the matrix $A$
appearing in \eqref{hbvm_rk1} (i.e., (\ref{rk}), by virtue of
Theorem~\ref{collhbvm}), one obtains:
\begin{equation}
\label{rk_leg1} \P^{-1} A \P = \P^{-1} \mathcal A \P \pmatrix{ll} I_s \\
& O\endpmatrix = \pmatrix{ll} \widetilde{X}_s \\ & O
\endpmatrix,
\end{equation}

\no where $\widetilde X_s$ is the matrix defined in (\ref{tXs}).
Relation \eqref{rk_leg1} and \cite[Theorem\,5.11 on page 86]{HW}
prove that method (\ref{hbvm_rk1}) (i.e., HBVM$(k,s)$) satisfies
$C(s)$ and $D(s-1)$ and, hence, its order is $2s$.

\section{Conclusions}\label{fine}
In this paper, we have shown that the recently introduced class of
HBVMs$(k,s)$, when recast as Runge-Kutta method, have the matrix
of the corresponding Butcher tableau with the same nonzero
eigenvalues which, in turn, coincides with those of the matrix of
the Butcher tableau of the Gauss method of order $2s$, for all
$k\ge s$ such that $B(2s)$ holds.

Moreover, HBVM$(k,s)$ defined at the Gaussian nodes
$\{\tau_1,\dots,\tau_k\}$ on the interval $[0,1]$ are closely
related to the Gauss method of order $2k$ which, as is well known,
is a symplectic method.

An alternative proof of the order of convergence of HBVMs is also
provided.

\end{document}